\documentclass[12pt]{amsart}

\usepackage{amssymb}

\usepackage{enumerate}

\usepackage[all]{xy}

\usepackage{graphicx}

\makeatletter
\@namedef{subjclassname@2010}{%
  \textup{2010} Mathematics Subject Classification}
\makeatother

\newtheorem{thm}{Theorem}[section]
\newtheorem{cor}[thm]{Corollary}

\theoremstyle{definition}

\newtheorem{rem}[thm]{Remark}

\numberwithin{equation}{section}

\begin{document}


\baselineskip=17pt



\title[Cutting description...]{Cutting description of trivial 1-cohomology}

\author{Andrzej Czarnecki}
\address{Jagiellonian University \\ Faculty of Mathematics and Computer Science \\ \L ojasiewicza 6,\, 30-348,\, Krakow Poland}
\email{andrzejczarnecki01@gmail.com}

\date{}

\begin{abstract}
A characterisation of trivial 1-cohomology for a broad class of metric spaces is presented. This trinket serves to remind that cohomology measures connectedness.
\end{abstract}



\maketitle

\section{Introduction.}
\label{In}

We will establish the following theorem

\begin{thm}\label{thm1} A connected and locally connected metric space $X$ has trivial first \v{C}ech cohomology group if and only if every connected open subset $U$ leaves $X\setminus U$ disconnected, provided it has a disconnected boundary.
\end{thm}

If we label the following conditions:

\begin{enumerate}
\item\label{1} $X$ is connected and locally connected;
\item\label{2} $H^1(X)=0$;
\item\label{3} $\partial U$ is disconnected;
\item\label{4} $X\setminus U$ is disconneced;
\end{enumerate}

then Theorem \ref{thm1} accounts for all nontrivial implications in 

\begin{thm}\label{thm2}
$(\ref{1})\Rightarrow\Big((\ref{2})\Leftrightarrow\Big(\forall U \text{ open and connected }  (\ref{3})\Leftrightarrow(\ref{4})\Big)\Big)$
\end{thm}

That $(\ref{4})$ always implies $(\ref{3})$, is an exercise on normality of metric spaces.


Of course, our theorems apply to the manifold category, and we can state one corollary in terms of de Rham cohomology, thus solving a PDE:

\begin{cor}
If every open domain $U$ of a manifold $M$ with $\partial U$ disconnected leaves $M\setminus U$ disconnected, then every equation

$$
df = \alpha
$$

has a solution, provided the 1-form $\alpha$ is closed.
\end{cor}

Throught this paper, $H^i(X)$ stands for the $i-$th reduced \v{C}ech cohomology group with constant $\mathbb{Z}$ coefficients.

This paper is a matured version of \cite{CKL} and, as mentioned there, observations of this kind have (minor) applications to complex analysis, concerning boundary of domains of holomorphy. Apart from the proof, we will give examples to show that the local connectedness cannot be ommited.

\section{The proof.}

All the algebraic topology material used here is clasic and can be found in any of popular textbooks on the subject.
Recall that the $0-$th group $H^0(A)$ is always free and in the locally connected setting its rank is equal to the number of connected components of $A$ minus 1. It is natural to apply the Mayer-Vietoris sequence in any problem concerning decompositions of a space and cohomology. However, some care in our case is needed. Recall that for every pair of open sets $A$ and $B$, covering the space $X$, we have an exact sequence

\begin{displaymath}
\xymatrix{ H^{0}(X) \ar[r] & H^{0}(A)\oplus H^{0}(B) \ar[r] & H^{0}(A\cap B) \ar[r]^{\partial_*} & H^{1}(X) \ar[r] &\ldots }
\end{displaymath}

We label only the so-called connecting homomorphism for future reference. To put ourselves in such setting, we consider small open neighbourhoods of $X\setminus U$, closure of $U$, $\overline{U}$ and boundary of $U$, $\partial U$: $X\setminus U_{\lambda}$, $\overline{U}_{\lambda}$, $\partial U_{\lambda}$, respectively.

\begin{displaymath}
\xymatrix{ H^{0}(X) \ar[r] & H^{0}(X\setminus \overline{U}_{\lambda})\oplus H^{0}(\overline{U}_{\lambda}) \ar[r] & H^{0}(\partial \overline{U}_{\lambda}) \ar[r]^{\partial_*} & H^{1}(X) \ar[r] & \ldots }
\end{displaymath}


The directed system of such neighbourhoods converges to our initial sets, and this is reflected by convergence in cohomology, by rigidity of \v{C}ech cohomology in metric spaces. Thus we have an exact sequence

\begin{displaymath}
\xymatrix{ H^{0}(X) \ar[r] & H^{0}(X\setminus U)\oplus H^{0}(\overline{U}) \ar[r] & H^{0}(\partial U) \ar[r]^{\partial_*} & H^{1}(X) \ar[r] & \ldots }
\end{displaymath}

Assume now that $X$ is connected, $U$ is a domain and $H^{1}(X)$ is trivial. The sequence takes form

\begin{displaymath}
\xymatrix{ 0 \ar[r] & H^{0}(X\setminus U) \ar[r] & H^{0}(\partial U) \ar[r] & 0 }
\end{displaymath}

This establishes the bijection between the components of the boundary and of the complement, and thus one implication in our theorem.

\begin{rem}
Dropping the assumption that $U$ is connected, we still get an exact sequence

\begin{displaymath}
\xymatrix{ H^{0}(X) \ar[r] & H^{0}(X\setminus U)\oplus H^{0}(U) \ar[r] & H^{0}(\partial U) \ar[r]^{\partial_*} & Im\partial_* \ar[r] & 0 }
\end{displaymath}

Note that $H^1(X)$, and thus also $Im\partial_*$, are free groups. Exactness means that the alternating sum of ranks of the groups in the sequence (its Euler characteristic) is zero

$$rk\,H^{0}(X) - rk\,H^{0}(X\setminus U) - rk\,H^{0}(U) + rk\,H^{0}(\partial U) + rk\,Im\partial_* = 0$$

Translating that to the number of connected components (we write $\#A$ for number of connected components of $A$), when $X$ is connected, we get inequalities

$$-1\leq \#\partial U - \#X\setminus U - \#U \leq-1 + rk\,H^1(X)$$

Note that for a broad class of spaces (spaces with ``good'' coverings in the sense of homotopy theory, manifolds, for example) $rk\,H^1(X)$ is bounded by $rk\,\pi_1(X)$.

\end{rem}

For the other implication, assume $H^{1}(X)$ nontrivial. We will find a domain with a connected complement and disconnected boundary.

$H^{1}(X) = \varinjlim H^{1}\mathcal{U}$, injective limit with respect to the directed system of all open coverings of $X$ -- without loss of generality, coverings by connected sets. Hence a nontrivial class in $H^{1}(X)$ arises as a nontrivial class in some $H^{1}(\mathcal{V})$ (and in all of it's refinements). $H^{1}(\mathcal{V})$ is in turn equal to $H_{S}^{1}(N\mathcal{V})$, singular cohomology of the nerve of $\mathcal{V}$, which is a simplicial complex. We can assume that $N\mathcal{V}$ is truncated over dimension 2, since we are interested only in the first cohomology group.
For any simplicial complex $K$, there is a 1:1 corespondence between $H_{S}^{1}(K)$ and $[K,\mathbb{S}^{1}]$, the homotopy classes of continuous maps from $K$ to the circle. Therefore a nontrivial class in $H_{S}^{1}(N\mathcal{V})$ is represented by a map $\theta$ from $N\mathcal{V}$ to $\mathbb{S}^1$. This map can be chosen simplicial (for a sufficiently fine simplicial structure on the cirlce; note that simplicial circle has at least three vertices) and without local extrema (a point $x$ is a local extremum of $\theta:N\mathcal{V}\rightarrow \mathbb{S}^1$ if it is a genuine local extremum in a neighbourhood $V_x$ of $\theta|_{V_x}\rightarrow B(\theta(x),\epsilon)\subset\mathbb{S}^1$, the small ball in $\mathbb{S}^1$ identified with an interval in $\mathbb{R}$).
Starting from any vertex, enumerate the vertices in the circle clockwise. Pick any vertex $v_n\in\mathbb{S}^1$. The vertices $a_i$ in $\theta^{-1}(v)$ are open sets in the covering $\mathcal{V}$. Any connected component of $\bigcup a_i$ must have disconnected boundary (disconnected by disjoint open sets $\theta^{-1}(v_{n-1})$ and $\theta^{-1}(v_{n+1})$). Moreover, there exists at least one connected component $A\subset\bigcup a_i$ such that its complement has a connected component meeting both $\theta^{-1}(v_{n-1})$ and $\theta^{-1}(v_{n+1})$ (otherwise $\theta$ would be nullhomotopic). For such $A$, pick this component of the complement, $B$. 
The domain $U=\bigcup\{v\text{ vertex in $N\mathcal{V}$}\,\,|\, v\cap B=\emptyset\}$ has a disconnected boundary and a connected complement. This proves the other implication in our theorem.

\section{Counterexamples.}

As for the counterexample concerning local connectedness, consider "rational Hawaiian earring", a dense subspace of a ball in $\mathbb{R}^2$:

$$\mathcal{H}_{\mathbb{Q}}=\bigcup \partial \mathbb{B}((0,q),q)$$

sum taken over all positive rational numbers up to 1. This connected space has obviously nontrivial 1-cohomology, and all its domains must contain the point $(0,0)$. The only case of such domain having a connected boundary is precisely when the boundary is equal to the complement and is contained in one of the circles. Thus, in terms of Theorem \ref{thm2}, $(\ref{3})\Leftrightarrow(\ref{4})$ does not imply $(\ref{2})$ without local connectedness.

We note however that trivial 1-cohomology always implies bijection between quasi-components of the complement and of the boundary of a domain, because the sequence

\begin{displaymath}
\xymatrix{ 0 \ar[r] & H^{0}(X\setminus U) \ar[r] & H^{0}(\partial U) \ar[r] & 0 }
\end{displaymath}

remains exact, and -- without the local connectedness assumed -- the ranks of the groups measure the number of quasi-components (minus 1). 

To finish, we note that the last remaining one-directional implication in Theorem \ref{thm2} cannot be reversed by the following counterexample.



The Knaster-Kuratowski fan (a cone over rationals) is contractible and satisfies the cutting condition (for reasons similar as in the case of Hawaiian earring), but is not locally connected. Observe, however, that this space is not locally homogenous (the vertex is topologically different from other points).

\subsection*{Acknowledgements}
Author thanks Marek Jarnicki for requesting a solution of this problem. Author is supported by IPhDPP of the FNP cofinanced by the EU under ERDF.

\end{document}